\documentclass{article}

\usepackage{amssymb}
\usepackage{amsmath}
\usepackage{amsfonts}
\usepackage{mathtools}
\usepackage{amsthm}
\newtheorem{theorem}{Theorem}

\newtheorem{corollary}[theorem]{Corollary}

\newtheorem{definition}[theorem]{Definition}

\newtheorem{lemma}[theorem]{Lemma}

\newtheorem{proposition}[theorem]{Proposition}
\newtheorem{remark}[theorem]{Remark}

\usepackage{authblk}

\title{Dichotomies uniform on subspaces\\ and formulas for dichotomy spectra}
\author[1]{Adam Czornik\footnote{adam.czornik@polsl.pl}}
\author[2]{Konrad Kitzing\footnote{konrad.kitzing@tu-dresden.de}}
\author[2]{Stefan Siegmund\footnote{stefan.siegmund@tu-dresden.de}}
\affil[1]{Faculty of Automatic Control, Electronics and Computer Science, Silesian University of Technology, Gliwice, Poland}
\affil[2]{Institute of Analysis, Faculty of Mathematics, TU Dresden, Germany}
\date{\today}

\newcommand{\N}{\mathbb{N}}
\newcommand{\Z}{\mathbb{Z}}
\newcommand{\R}{\mathbb{R}}
\newcommand{\T}{\mathbb{T}}
\allowdisplaybreaks
\usepackage{enumitem}
\setlist[enumerate]{label*=(\alph*),ref=(\alph*)}

\begin{document}

\maketitle

\begin{abstract}
We introduce a notion of dichotomy which generalizes the classical concept of exponential dichotomy and the recent notion of Bohl dichotomy. A key attribute is the discussion of the sets of subspaces of the state space on which the dichotomy estimates  are uniform. Two main results are a dichotomy spectral theorem based on our notion of dichotomy which is uniform on subspaces and a formula for the dichotomy spectral intervals which is new for the Bohl dichotomy spectrum as well as for the classical exponential dichotomy spectrum.

\end{abstract}


\section{Introduction}\label{sec:intro}

In this paper we develop formulas for the dichotomy spectrum of linear time-varying difference equations on $\mathbb{R}^d$,
\(d \in \N_{>0}\)
\begin{equation}\label{1}
    x(n+1) = A(n)x(n),
    \qquad
    n \in \T,
\end{equation}
with one-sided time \(\T = \N\) or two-sided time \(\T = \Z\) and invertible coefficients \(A \colon \T \to \mathrm{GL}(\R^d)\) such that \(A\) and \(A^{-1} \colon\T \to \mathrm{GL}(\R^d)\), \(n \mapsto A(n)^{-1}\) are bounded.

There is a rich body of literature on hyperbolicity concepts for systems \eqref{1}, most prominently exponential dichotomy, but also more recently Bohl dichotomy, see Definition \ref{ED-BD} and \cite{Czornik2023} (for the case $\T=\N$) and Remark \ref{rem:dich-history} for a short history and references. In this introduction we observe that a main difference between these known dichotomy notions is the set of subspaces on which the dichotomy estimates for the solutions of \eqref{1} are uniform. This will lead us in Definition \ref{def:dicho} to a new notion of dichotomy which is uniform on subspaces and which generalizes the concepts of exponential and Bohl dichotomy.

Section 2 contains results on properties of dichotomous systems \eqref{1} as a preparation for the main results in Section 4. In Section 3 we recall in Definition \ref{def:BohlExp} Bohl exponents from \cite{Czornik2023} and investigate their relation to our general notion of dichotomy.
Section 4 contains the two main results, a dichotomy spectral theorem based on our general notion of dichotomy which is uniform on subspaces (Theorem \ref{thm:spectral}) and a formula for the dichotomy spectral intervals (Theorem \ref{thm:spectral_formula}). 

For an initial value \(x_0 \in \R^d\) the solution $x(\cdot,x_0) : \T \to \R^d$ of \eqref{1} starting at $0$ in $x_0$ satisfies \(x(0,x_0) = x_0\),
\begin{equation*}
   x(n,x_0) = A(n-1) \cdots A(0)x_0,
   \qquad n > 0,
\end{equation*}
and if $\T = \Z$
\begin{equation*}
   x(n,x_0) = A(n)^{-1} \cdots A(-1)^{-1}x_0,
   \qquad n < 0.
\end{equation*}
Dichotomies are solution estimates on subspaces of the state space $\R^d$ which form a splitting. Denote for a subspace \(L \subseteq \R^d\) the set of subspaces of $L$ by
\begin{align*}
    \mathcal G(L) &\coloneqq \{U \subseteq L : \text{\(U\) subspace}\}.
\end{align*}
For \(L_1,L_2 \in \mathcal G(\mathbb{R}^d)\) we say that \((L_1,L_2)\) is a splitting of \(\R^d\) if  \(L_1\oplus L_2 = \R^d\).

Let \((L_1, L_2)\) be a splitting of \(\R^d\).
We denote by \(\pi_{L_1}\) the projection onto \(L_1\) along \(L_2\). Then \(\pi_{L_2} \coloneqq I - \pi_{L_1}\) is the projection onto \(L_2\) along \(L_1\), where \(I\) denotes the identity mapping on \(\R^d\).
For \(V \in \mathcal{G}(\R^{d})\) and \(i \in \{1,2\}\)
\begin{equation*}
    \pi_{L_i}[V] \coloneqq \{\pi_{L_i} v : v \in V\}
\end{equation*}
is the projection of \(V\) onto \(L_i\). For the following definition see e.g.\ \cite{Czornik2023} (or \cite{Coppel1978}; note that system \eqref{1} has bounded growth) and the references therein.

\begin{definition}[Exponential and Bohl dichotomy]\label{ED-BD}
Let \((L_1,L_2)\) be a splitting of \(\R^d\). System \eqref{1} has a \emph{Bohl dichotomy on $(L_1, L_2)$} if there exists an $\alpha > 0$ and functions $C_1, C_2 : \mathbb{R}^d \to (0,\infty)$ such that
\begin{align*}
   \|x(n,x_0)\|
   &\leq
   C_1(x_0) e^{-\alpha(n-m)} \|x(m,x_0)\|,
   &\quad m, n \in \T, n \geq m, x_0 \in L_1,
\\
   \|x(n,x_0)\|
   &\geq
   C_2(x_0) e^{\alpha(n-m)} \|x(m,x_0)\|,
   &\quad m, n \in \T, n \geq m,x_0 \in L_2.
\end{align*}
In case $C_1$ and $C_2$ can be chosen independent of $x_0$, system \eqref{1} has an \emph{exponential dichotomy}.
\end{definition}

\begin{remark}[History of dichotomy notion]
\label{rem:dich-history}\hfill

(a) Foundations.
Massera and Schaeffer \cite{Massera1,Massera2} coined the name \emph{exponential dichotomy} and formulated the canonical definition of exponentially dichotomous linear differential equations. The foundational formulation of the theory of exponentially dichotomous systems was completed through the monographs of Massera and Schaeffer \cite%
{Massera1966} and Daleckii and Krein \cite{DaleckiKrein1974}, which summarized earlier results on exponentially dichotomous systems and outlined paths for the
development of this theory by presenting numerous open questions.

For difference equations, this concept first appeared in the work of Coffman
and Schaeffer \cite{Coffman Schaffer 1967}, where the notion of uniform
exponential dichotomy for infinite-dimensional, one-sided discrete-time systems was introduced, highlighting differences between continuous and discrete-time
systems due to the potential non-invertibility of the transition operator.
The first definitions of dichotomy for two-sided sequences were provided in  \cite{Boor 1980} and \cite[Definition 7.6.4]{Henry 1981}. Subsequently,
many authors have studied the concept of uniform exponential dichotomies
for discrete systems with both invertible and non-invertible coefficients
\cite{Ben-Artzi Gohberg 1989, Ben-Artzi Gohberg 1990, Aulbah Minh Zabreiko 1994, Halanay
Ionescu 1994}.

The origins of the concept go back to Perron’s seminal work \cite{Perron1930}, one of the earliest studies of this property. For comprehensive treatments of exponential dichotomy and its applications, see \cite{Coppel1978, Lin2000, Mitropolsky2002}.

(b) Concepts of uniformity.
The term \emph{uniform} usually refers to the fact that the dichotomy estimates are uniform with respect to both time and the initial condition. While the role of uniformity with respect to time has been extensively analyzed in the
literature (see e.g.\ \cite{Barreira Pesin 2006, Barreira 2008, Barreira2018, Silva 2021} and references therein),  research of the significance of uniformity with respect to the initial condition has only
recently been initiated (see e.g.\ \cite{Bekryaeva2010, Doan2017, Barreira2018b, Czornik2023, Czornik2024}). These studies examine a type of dichotomy, referred to as \emph{Bohl
dichotomy} in \cite{Czornik2023}, where the convergence with
respect to time is uniform, but uniformity with respect to the initial condition is restricted to each one-dimensional subspace.

(c) Weak dichotomies.
The version of exponential dichotomy where the constants $C_1$, $C_2$ in Definition \ref{ED-BD} may depend on the initial condition was first introduced for continuous-time systems in \cite{Bekryaeva2010}, without a specific name. It was later explicitly termed weak exponential dichotomy in \cite{Barabanov2020}. For discrete-time systems on the half-line, the same notion appeared in \cite{Barreira2018b} and has since been systematically studied under the name Bohl dichotomy in \cite{Czornik2023, Czornik2024}.
\end{remark}

In the following remark we emphasize the subspaces on which the dichotomy estimates in Definition \ref{ED-BD} hold uniformly. Denote for \(k \in \N\) and a subspace \(L \subseteq \R^d\) the set of $k$-dimensional subspaces of $L$ by
\begin{align*}
    \mathcal G_{k}(L) &\coloneqq \{U \in \mathcal G(L) : \dim U = k\}.
\end{align*}

\begin{remark}[Uniformity subspaces of Bohl and exponential dichotomy]\label{rem:uni_sub_Bohl_expo}
If \eqref{1} has a Bohl dichotomy on $(L_1,L_2)$ then, using the homogeneity of the solution $x(\cdot, c x_0) = c x(\cdot, x_0)$, $c \in \R$, the dichotomy estimates hold uniformly on all one-dimensional subspaces of $L_1$ and $L_2$, respectively. More precisely, for all \(U \in \mathcal G_1(L_1)\), there exists \(C(U) > 0\) with
\begin{equation*}
    \Vert x(n,x_0)\Vert \leq C(U)\mathrm e^{-\alpha(n-m)}\Vert x(m,x_0)\Vert,
    \qquad
    m,n\in\T,\, n\geq m,\, x_0 \in U,
\end{equation*}
and for all \(U \in \mathcal G_1(L_2)\), there exists \(C(U) > 0\) with
\begin{equation*}
    \Vert x(n,x_0)\Vert \geq C(U)\mathrm e^{\alpha(n-m)}\Vert x(m,x_0)\Vert,
    \qquad
    m,n\in\T,\, n\geq m,\, x_0 \in U.
\end{equation*}
If \eqref{1} has an exponential dichotomy on $(L_1,L_2)$ then the dichotomy estimates hold uniformly on all subspaces of $L_1$ and $L_2$, respectively. \end{remark}

We introduce a new dichotomy notion based on dichotomy estimates which hold uniformly on a set of subspaces of the splitting which covers the splitting.

\begin{definition}[Dichotomy uniform on subspaces of a splitting]\label{def:dicho}
Let \((L_1,L_2)\) be a splitting of \(\R^d\) and \(\mathcal U_1 \subseteq  \mathcal G(L_1)\), \(\mathcal U_2 \subseteq  \mathcal G(L_2)\) be covers of $L_1$ and $L_2$, respectively, i.e.
\begin{equation}\label{eq:splitting_admissble}
   \bigcup \mathcal U_1 = L_1
   \qquad \text{and} \qquad
   \bigcup \mathcal U_2 = L_2.
\end{equation}
We say that system \eqref{1} has a \emph{dichotomy on \((L_1,L_2)\) uniformly on the subspaces in \((\mathcal U_1,\mathcal U_2)\)} if there exists an $\alpha > 0$ such that for each \(U_1 \in \mathcal U_1\), there exists \(C(U_1) > 0\) with
\begin{equation}\label{def:dicho_B_1}
    \Vert x(n,x_0)\Vert \leq C(U_1)\mathrm e^{-\alpha(n-m)}\Vert x(m,x_0)\Vert,
    \qquad
    m,n\in\T,\, n\geq m,\, x_0 \in U_1,
\end{equation}
and for each \(U_2 \in \mathcal U_2\), there exists \(C(U_2) > 0\) with
\begin{equation}\label{def:dicho_B_2}
    \Vert x(n,x_0)\Vert \geq C(U_2)\mathrm e^{\alpha(n-m)}\Vert x(m,x_0)\Vert,
    \qquad
    m,n\in\T,\, n\geq m,\, x_0 \in U_2.
\end{equation}
We say that system \eqref{1} has a \emph{dichotomy on \((L_1,L_2)\)} if there exist covers \(\mathcal U_1 \subseteq  \mathcal G(L_1)\) of $L_1$ and \(\mathcal U_2 \subseteq  \mathcal G(L_2)\) of $L_2$ such that system \eqref{1} has a dichotomy on \((L_1,L_2)\) uniformly on the subspaces in \((\mathcal U_1,\mathcal U_2)\).
\end{definition}

\begin{remark}[Examples of dichotomies uniform on subspaces of a splitting]\label{rem:dich_exmpl}
By Definition \ref{ED-BD} Bohl and exponential dichotomy are examples of dichotomies uniform on subspaces. If system \eqref{1} has a dichotomy uniform on subspaces \((\mathcal U_1,\mathcal U_2)\) of a splitting \((L_1,L_2)\) then it has an

(a) exponential dichotomy if $\mathcal U_1 = \{L_1\}$ and $\mathcal U_2 = \{L_2\}$ and a

(b) Bohl dichotomy if $\mathcal U_1 = \mathcal G_1(L_1)$ and $\mathcal U_2 =  \mathcal G_1(L_2)$.

\end{remark}

The important and obvious observation that a dichotomy estimate which holds uniformly on a subspace $U$ also holds uniformly on all subspaces of $U$ is formulated in the following theorem.

\begin{theorem}[Refining uniformity subspaces of dichotomy]
\label{thm:dicho-refine}
If system \eqref{1} has a dichotomy on \((L_1,L_2)\) uniformly on the subspaces in \((\mathcal U_1,\mathcal U_2)\), then it also has a dichotomy on \((L_1,L_2)\) uniformly on the subspaces in the refinement \((\mathcal V_1,\mathcal V_2)\) with
\begin{equation*}
   \mathcal V_1
   \coloneqq
   \bigcup_{U \in \mathcal U_1}
   \mathcal G(U)
   \qquad \text{and} \qquad
   \mathcal V_2
   \coloneqq
   \bigcup_{U \in \mathcal U_2}
   \mathcal G(U).
\end{equation*}
\end{theorem}

\begin{proof}
Let $i \in \{1,2\}$.
The set $\mathcal V_i$ consists of all subspaces of the subspaces in $\mathcal U_i$ and therefore it also covers $L_i$. Let $V_i \in \mathcal V_i$. Then there exists $U_i \in \mathcal U_i$ such that $V_i$ is a subspace of $U_i$ and the dichotomy estimate which holds on $U_i$ also holds on $V_i$.
\end{proof}

In Remark \ref{rem:uni_sub_Bohl_expo} it was pointed out that for Bohl and exponential dichotomies the estimates are uniform on all one-dimensional subspaces of the corresponding spaces of the splitting. This property also holds for the new dichotomy notion defined in Definition \ref{def:dicho}.

\begin{remark}[Dichotomy is uniform on one-dimensional subspaces]\label{rem:refinement}
If system \eqref{1} has a dichotomy on \((L_1,L_2)\) uniformly on the subspaces in \((\mathcal U_1,\mathcal U_2)\), then it also has a dichotomy on \((L_1,L_2)\) uniformly on the subspaces in \((\mathcal G_1(L_1),\mathcal G_1(L_2))\). This is due to the fact that if \(L \neq \{0\}\) is a subspace of \(\R^d\) and \(\mathcal U \subseteq  \mathcal G(L)\) is a cover of $L$ then
\begin{equation*}
   \mathcal G_1(L) \subseteq  \bigcup_{U \in \mathcal U}
   \mathcal G(U).
\end{equation*}
As a consequence of Remark \ref{rem:dich_exmpl}, if system \eqref{1} has a dichotomy on \((L_1,L_2)\) then it also has a Bohl dichotomy on \((L_1,L_2)\).
\end{remark}

The notions of exponential and Bohl dichotomy spectrum (see e.g.\ \cite{Czornik2023}) are based on dichotomies of the \emph{$\gamma$-shifted} systems for $\gamma \in \mathbb{R}$
\begin{equation}\label{gamma-shift}
   x(n+1)=e^{-\gamma }A(n)x(n)
   , \qquad
   n \in \mathbb{T}.
\end{equation}

\begin{theorem}[Exponential and Bohl dichotomy spectrum]
\label{thm:exp-bohl-spectrum}
The Bohl dichotomy spectrum of \eqref{1}
\begin{equation*}
   \Sigma _{\mathrm{BD}}(A)
   \coloneqq
   \left\{ \gamma \in \mathbb{R}:\eqref{gamma-shift}%
    \text{ has no Bohl dichotomy}\right\},
\end{equation*}
as well as the exponential dichotomy spectrum of \eqref{1}
\begin{equation*}
   \Sigma _{\mathrm{ED}}(A)
   \coloneqq
   \left\{ \gamma \in \mathbb{R}:\eqref{gamma-shift}%
    \text{ has no exponential dichotomy}\right\},
\end{equation*}
is the nonempty union of at most \(d\) compact intervals.
\end{theorem}

In this paper we generalize Theorem \ref{thm:exp-bohl-spectrum} to dichotomies which hold uniformly on subspaces of prescribed dimensions and provide formulas for the endpoints of the spectral intervals (Theorems \ref{thm:spectral} and \ref{thm:spectral_formula}) which are also new for the exponential and Bohl dichotomy spectrum.

\section{Dichotomies}

In this section we study dichotomies which are uniform on subspaces of a splitting. The estimates \eqref{def:dicho_B_1} and \eqref{def:dicho_B_2} are the building blocks of Definition \ref{def:dicho} which we abbreviate for convenience in the following definition.

\begin{definition}[Dichotomy estimates]
\label{def:D}
Let \(\gamma \in \R\) and \(U \in \mathcal G(\R^d)\).
\begin{enumerate}
\item
We say \(\mathrm D_1(\gamma,U)\) holds if there exists \(C \colon U \to \R_{>0}\), such that
\begin{equation*}
    \Vert x(n,x_0)\Vert \leq C(x_0)\mathrm e^{\gamma(n-m)}\Vert x(m,x_0)\Vert,
    \qquad
    m,n\in\T,\, n\geq m,\, x_0 \in U.
\end{equation*}
If the mapping \(C\) is constant, we say that \(\mathrm D_1(\gamma,U)\) holds uniformly.
\item
We say \(\mathrm D_2(\gamma,U)\) holds if there exists \(C \colon U \to \R_{>0}\), such that
\begin{equation*}
    \Vert x(n,x_0)\Vert \geq C(x_0)\mathrm e^{\gamma(n-m)}\Vert x(m,x_0)\Vert,
    \qquad
    m,n\in\T,\, n\geq m,\, x_0 \in U.
\end{equation*}
If the mapping \(C\) is constant, we say that \(\mathrm D_2(\gamma,U)\) holds uniformly.
\end{enumerate}
\end{definition}

Using the notation of Definition \ref{def:D} system \eqref{1} has a dichotomy with splitting \((L_1,L_2)\), if and only if there exists \(\alpha > 0\) and a covering \(\mathcal U_1\) of \(L_1\) and a covering \(\mathcal U_2\) of \(L_2\), such that for all \(U_1 \in \mathcal U_1\) and \(U_2 \in \mathcal U_2\) the properties \(\mathrm D_1(-\alpha,U_1)\) and \(\mathrm D_2(\alpha,U_2)\) hold uniformly. In particular, in this case $C$ can be chosen as a constant function.

\begin{remark}[Characterization of dichotomy of $\gamma$-shifted system]
For \(\gamma \in \R\), the solution of the \(\gamma\)-shifted system \eqref{gamma-shift} which starts at $0$ in \(x_0 \in \R^d\) is
\(\big(\mathrm e^{-\gamma n}x(n,x_0)\big)_{n\in\T}\).
The \(\gamma\)-shifted system has a dichotomy with splitting \((L_1,L_2)\) if and only if there exists \(\alpha > 0\) and a covering \(\mathcal U_1\) of \(L_1\) and a covering \(\mathcal U_2\) of \(L_2\), such that \(\mathrm D_1(\gamma-\alpha,U)\) holds uniformly for all \(U \in \mathcal U_1\) and \(\mathrm D_2(\gamma+\alpha,U)\) holds uniformly for all \(U \in \mathcal U_2\).
\end{remark}

\begin{proposition}[Dichotomy estimates are dichotomous]\label{prop:est_dicho}
Let \(\alpha > 0\), \(\gamma \in \R\) and \(U_1,U_2 \in \mathcal G(\R^d)\).
If \(\mathrm D_1(\gamma-\alpha,U_1)\) and \(\mathrm D_2(\gamma+\alpha,U_2)\) hold, then \(U_1 \cap U_2 = \{0\}\).

More generally, suppose that \(U \in \mathcal G(\R^d)\), \(\alpha>0\) and \(\gamma \in \R\).
If \(\mathrm D_2(\gamma+\alpha,U)\) holds and if \(\big(\mathrm e^{-\gamma n}x(n,x_0)\big)_{n\in\N}\) is bounded for all \(x_0 \in U\), then \(U = \{0\}\).
\end{proposition}

\begin{proof}
Let \(x_0 \in U_1 \cap U_2\).
From \(\mathrm D_1(\gamma-\alpha,U_1)\) (with \(m = 0\) and noting that \(\N \subseteq \T\)), we conclude that there exists \(C_1(x_0)>0\) with
\begin{equation*}
    \Vert x(n,x_0)\Vert \leq C_1(x_0)\mathrm e^{(\gamma-\alpha)n}\Vert x_0\Vert,\qquad n \in \N.
\end{equation*}
In particular, \(\big(\mathrm e^{-\gamma n}x(n,x_0)\big)_{n\in\N}\) is bounded.
The estimate \(\mathrm D_2(\gamma + \alpha,U_2)\) yields a \(C_2(x_0)>0\) with
\begin{equation*}
    \Vert x_0\Vert \leq C_2(x_0)\mathrm e^{-(\gamma+\alpha)n}\Vert x(n,x_0)\Vert,\qquad n \in \N.
\end{equation*}
Letting \(n \to \infty\), using the fact that \(\alpha>0\) and the boundedness of \(\big(\mathrm e^{-\gamma n}x(n,x_0)\big)_{n\in\N}\), we conclude \(x_0 = 0\).
\end{proof}

\begin{proposition}[Monotonicity of dichotomy subspaces]\label{prop:mono_dicho}
Let \(\gamma,\widetilde\gamma \in \R\) with \(\gamma \leq \widetilde\gamma\).
Suppose that the \(\gamma\)-shifted system \eqref{gamma-shift} has a dichotomy with splitting \((L_1,L_2)\) and that the \(\widetilde\gamma\)-shifted system has a dichotomy with splitting \((\widetilde L_1, \widetilde L_2)\).
Then \(L_1 \subseteq \widetilde L_1\).
If \(\T = \Z\), then \(L_2 \supseteq \widetilde L_2\).
\end{proposition}

\begin{proof}
We prove \(L_1 \subseteq \widetilde L_1\).
Let \(x_1 \in L_1\).
There are \(\widetilde x_1 \in \widetilde L_1\), \(\widetilde x_2 \in \widetilde L_2\) with \(x_1 = \widetilde x_1 + \widetilde x_2\).
We show \(\widetilde x_2 = 0\) by applying Proposition \ref{prop:est_dicho} to the \(\widetilde\gamma\)-shifted system, so that \(x_1 = \widetilde x_1\), i.e.\ \(L_1 \subseteq \widetilde L_1\).
Indeed, there are \(\alpha_1,\alpha_2 > 0\), such that \(D_1(\gamma-\alpha_1,\operatorname{span}\{x_1\})\) and \(D_1(\widetilde\gamma-\alpha_2,\operatorname{span}\{\widetilde x_1\})\) and \(\mathrm D_2(\widetilde\gamma+\alpha_2,\operatorname{span}\{\widetilde x_2\})\) hold by Remark \ref{rem:refinement}.
In particular, for \(\alpha \coloneqq \min(\alpha_1,\alpha_2)\) there is \(C>0\), such that for \(n \in \N\),
\begin{align*}
    \Vert \mathrm e^{-\widetilde\gamma n}x(n,\widetilde x_2)\Vert
    &=
    \mathrm e^{-\widetilde\gamma n}\Vert x(n,x_1 - \widetilde x_1)\Vert
    \\
    &\leq
    \mathrm e^{-\widetilde\gamma n}\big(C\mathrm e^{(\gamma-\alpha)(n-0)}\Vert x(0,x_1)\Vert
    +
    C\mathrm e^{(\widetilde\gamma-\alpha)(n-0)}\Vert x(0,\widetilde x_1)\Vert\big)
    \\
    &=
    C\big(\mathrm e^{(\gamma - \widetilde\gamma - \alpha)n}\Vert x_1\Vert + \mathrm e^{-\alpha n}\Vert \widetilde x_1\Vert\big).
\end{align*}
Since \(\gamma - \widetilde\gamma \leq 0\), we see that \(\big(\mathrm e^{-\widetilde\gamma n}x(n,\widetilde x_2)\big)_{n\in\N}\) is bounded, so that Proposition \ref{prop:est_dicho} yields \(\operatorname{span}\{\widetilde x_2\} = \{0\}\), i.e.\ \(\widetilde x_2 = 0\).
That \(L_2 \supseteq \widetilde L_2\) if \(\T = \Z\), can be shown similarly.
\end{proof}

\begin{proposition}[Uniqueness and dynamic characterization of dichotomy subspaces]\label{prop:uni_rep_dicho_spaces}
Let \(\gamma \in \R\).
Suppose that the \(\gamma\)-shifted system \eqref{gamma-shift} has a dichotomy with splitting \((L_1,L_2)\).
Then
\begin{equation*}
    L_1 = \{x_0 \in \R^d : \lim_{n\to\infty}\mathrm e^{-\gamma n}x(n,x_0) = 0\},
\end{equation*}
and if \(\T = \Z\), then
\begin{equation*}
    L_2 = \big\{x_0 \in \R^d : \lim_{n\to-\infty}\mathrm e^{-\gamma n}x(n,x_0) = 0\big\}.
\end{equation*}
Moreover, if for \(\widetilde\gamma \in \R\), the \(\widetilde\gamma\)-shifted system has a dichotomy with splitting \((\widetilde L_1,\widetilde L_2)\) and \(\dim L_1 = \dim \widetilde L_1\), then \(L_1 = \widetilde L_1\) and if \(\T = \Z\), then \(L_2 = \widetilde L_2\).
\end{proposition}

\begin{proof}
We prove the representation of \(L_1\) only and set \(M \coloneqq \{x_0 \in \R^d : \lim_{n\to\infty}\mathrm e^{-\gamma n}x(n,x_0) = 0\}\).
\(L_1 \subseteq M\), since there is \(\alpha > 0\), such that \(D_1(\gamma-\alpha,\operatorname{span}\{x_0\})\) holds for all \(x_0 \in L_1\).
To show \(M \subseteq L_1\), let \(x_0 \in M\) and let \(x_1 \in L_1\), \(x_2 \in L_2\) with \(x_0 = x_1 + x_2\).
Similar to the proof of \(\widetilde x_2 = 0\) in the proof of Proposition \ref{prop:mono_dicho}, it follows that \(x_2 = 0\) and thus \(x_0 = x_1 \in L_1\).
Now suppose \(\widetilde\gamma \in \R\) and assume w.l.o.g.\ \(\widetilde\gamma \leq \gamma\).
From Proposition \ref{prop:mono_dicho}, we obtain that \(\widetilde L_1 \subseteq L_1\) and since \(\dim \widetilde L_1 = \dim L_1\), also \(\widetilde L_1 = L_1\).
Similarly one can show that \(\widetilde L_2 = L_2\) if \(\T = \Z\).
\end{proof}

\section{Bohl exponents and dichotomy estimates}

A useful tool for studying properties of dichotomies is the concept of Bohl exponents which were introduced for individual solutions in \cite{Bohl1913} and later extended to subspace exponents in \cite{Ben-Artz Gohberg 1991} and \cite{Czornik2023}.
In this section we recall the upper and lower Bohl exponents from \cite{Czornik2023} and investigate their relation to the dichotomy estimates in Definition \ref{def:D}.
We set \(\sup\emptyset \coloneqq -\infty\) and \(\inf\emptyset \coloneqq \infty\).

\begin{definition}[Upper and lower Bohl exponent]\label{def:BohlExp}
For \(U \in \mathcal G(\R^d)\) we define the \emph{upper Bohl exponent}
\begin{equation*}
    \overline\beta(U) \coloneqq \inf_{N\in\N}\sup\bigg\{\frac 1{n-m}\ln\frac{\Vert x(n,x_0)\Vert}{\Vert x(m,x_0)\Vert} : m,n\in\T,\,n-m>N,\,x_0\in U\setminus\{0\}\bigg\}
\end{equation*}
and the \emph{lower Bohl exponent}
\begin{equation*}
    \underline\beta(U) \coloneqq \sup_{N\in\N}\inf\bigg\{\frac 1{n-m}\ln\frac{\Vert x(n,x_0)\Vert}{\Vert x(m,x_0)\Vert} : m,n\in\T,\,n-m>N,\,x_0\in U\setminus\{0\}\bigg\}.
\end{equation*}
\end{definition}

It should be noted that in the literature a series of other but equivalent definitions of these
exponents can be found and sometimes they appear under different names (see \cite[Remark 8]{Czornik2023}).

\begin{remark}[Basic properties of Bohl exponents]
\label{rem:infty_Bohl_exponents}\hfill

(a) $\overline\beta(\{0\}) = -\infty$ and $\underline\beta(\{0\}) = \infty$.

(b) If \(U \neq \{0\}\) then  \(-\ln\Vert A^{-1}\Vert_\infty \leq \underline\beta(U) \leq \overline\beta(U) \leq \ln\Vert A\Vert_\infty\). In particular,
\begin{equation*}
   \underline\beta(U),
   \overline\beta(U) \in \R,
   \qquad
   U \in \mathcal G(\R^d)\setminus\{0\}.
\end{equation*}

(c) For \(U,V \in \mathcal G(\R^d)\)
\begin{align*}
   U \subseteq V
   \quad \Rightarrow \quad
   \underline\beta(U) \geq \underline\beta(V)
   \text{ and }
   \overline\beta(U) \leq \overline\beta(V).
\end{align*}
\end{remark}

The next two propositions state the relation between the upper and lower Bohl exponents and the dichotomy estimates \(D_1(\gamma,U)\) and \(D_2(\gamma,U)\), respectively.

\begin{proposition}[Upper Bohl exponent and dichotomy estimate \(\mathrm D_1(\gamma,U)\)]\label{prop:Rel_upper}
Let \(\mathcal V \subseteq \mathcal G(\R^d)\) be non-empty and for each \(V \in \mathcal V\) let \(\mathcal U(V) \subseteq \mathcal G(\R^d)\) be an arbitrary non-empty set of subspaces of $\R^d$.
Let \(\gamma \in \R\).
\begin{enumerate}
\item\label{prop:Rel_upper_1}
If there is \(V \in \mathcal V\), such that for all \(U \in \mathcal U(V)\) the estimate \(\mathrm D_1(\gamma,U)\) holds uniformly, then \(\gamma \geq \inf_{V\in\mathcal V}\sup_{U \in \mathcal U(V)}\overline\beta(U)\).
\item\label{prop:Rel_upper_2}
If \(\gamma > \inf_{V\in\mathcal V}\sup_{U \in \mathcal U(V)}\overline\beta(U)\), then there is \(V \in \mathcal V\), such that for all \(U \in \mathcal U(V)\) the estimate \(\mathrm D_1(\gamma,U)\) holds uniformly.
\end{enumerate}
\end{proposition}

\begin{proof}
We prove \ref{prop:Rel_upper_1}.
Let \(V \in \mathcal V\), such that \(\mathrm D_1(\gamma,U)\) holds for all \(U \in \mathcal U(V)\).
First consider \(U \in \mathcal U(V)\setminus\{0\}\).
We obtain by \(\mathrm D_1(\gamma,U)\) a constant \(C>0\), such that
\begin{equation*}
    \frac 1{n-m}\ln\frac{\Vert x(n,x_0)\Vert}{\Vert x(m,x_0)\Vert} \leq \frac 1{n-m}C + \gamma,
    \qquad
    m,n\in\T,\,n>m,x_0\in U\setminus\{0\}.
\end{equation*}
Hence for all \(N \in \N\), we obtain
\begin{equation*}
    \frac 1{n-m}\ln\frac{\Vert x(n,x_0)\Vert}{\Vert x(m,x_0)\Vert} \leq \frac CN + \gamma,
    \qquad
    m,n\in\T,\,n-m>N,x_0\in U\setminus\{0\}.
\end{equation*}
Thus \(\overline\beta(U)\leq\gamma\).
Also \(\overline\beta(\{0\}) = -\infty < \gamma\) and we obtain \(\sup_{U\in\mathcal U(V)}\overline\beta(U)\leq\gamma\).

We prove \ref{prop:Rel_upper_2}.
From \(\gamma > \inf_{V\in\mathcal V}\sup_{U \in \mathcal U(V)}\overline\beta(U)\) we conclude that there exists \(V \in \mathcal V\), such that for all \(U \in \mathcal U(V)\) we have
\begin{equation*}
    \gamma>\overline\beta(U) = \inf_{N\in\N}\sup\bigg\{\frac 1{n-m}\ln\frac{\Vert x(n,x_0)\Vert}{\Vert x(m,x_0)\Vert} : m,n\in\T,\,n-m>N,\,x_0\in U\setminus\{0\}\bigg\}.
\end{equation*}
Note that \(\mathrm D_1(\gamma,\{0\})\) always holds uniformly and for \(U \in \mathcal U(V)\setminus\{0\}\), there is \(N \in \N\) with
\begin{equation*}
    \gamma > \frac 1{n-m}\ln\frac{\Vert x(n,x_0)\Vert}{\Vert x(m,x_0)\Vert},
    \qquad
    m,n\in\T,\,n-m>N,\,x_0\in U\setminus\{0\}.
\end{equation*}
Rearranging yields \(\Vert x(n,x_0)\Vert \leq \mathrm e^{\gamma(n-m)}\Vert x(m,x_0)\Vert\) for \(m,n\in\T\), \(n-m>N\), \(x_0\in U\).
To conclude, we have to show that there is \(C>0\), such that
\begin{equation*}
    \Vert x(n,x_0)\Vert \leq C\mathrm e^{\gamma(n-m)}\Vert x(m,x_0)\Vert,
    \qquad
    m,n\in\T,\,n-m\leq N,\,x_0\in U.
\end{equation*}
Indeed, for \(x_0 \in U\setminus\{0\}\) this follows from
\begin{align*}
    \mathrm e^{-\gamma(n-m)}\frac{\Vert x(n,x_0)\Vert}{\Vert x(m,x_0)\Vert}
    &=
    \mathrm e^{-\gamma(n-m)}\frac{\Vert A(n-1)\cdots A(m)A(m-1)\cdots A(0)x_0\Vert}{\Vert A(m-1)\cdots A(0)x_0\Vert}
    \\
    &\leq
    \mathrm e^{-\gamma(n-m)}\Vert A(n-1)\cdots A(m)\Vert
\end{align*}
and by noting that \(\mathrm e^{-\gamma(n-m)}\Vert A(n-1)\cdots A(m)\Vert\) is bounded on \(\{(m,n)\in\T^2 : 0\leq n-m\leq N\}\), since \(A\) is bounded. Hence the inequality holds with $C
  :=
  \|A\|_\infty^N$.
\end{proof}

\begin{proposition}[Lower Bohl exponent and \(\mathrm D_2(\gamma,U)\)]\label{prop:Rel_lower}
Let \(\mathcal V \subseteq \mathcal G(\R^d)\) be non-empty and for every \(V \in \mathcal V\) let \(\mathcal U(V) \subseteq  \mathcal G(\R^d)\) be an arbitrary non-empty set of subspaces of $\R^d$.
Let \(\gamma \in \R\).
\begin{enumerate}
\item\label{prop:Rel_lower_1}
If there is \(V \in \mathcal V\), such that for all \(U \in \mathcal U(V)\) the estimate \(\mathrm D_2(\gamma,U)\) holds uniformly, then \(\gamma \leq \sup_{V\in\mathcal V}\inf_{U \in \mathcal U(V)}\underline\beta(U)\).
\item\label{prop:Rel_lower_2}
If \(\gamma < \sup_{V\in\mathcal V}\inf_{U \in \mathcal U(V)}\underline\beta(U)\), then there is \(V \in \mathcal V\), such that for all \(U \in \mathcal U(V)\) the estimate \(\mathrm D_2(\gamma,U)\) holds uniformly.
\end{enumerate}
\end{proposition}

\begin{proof}
The proof is along the lines of Proposition \ref{prop:Rel_upper}.
\end{proof}

\section{The spectral theorem}

In this section we generalize Theorem \ref{thm:exp-bohl-spectrum} and provide formulas for the spectrum. For each $\gamma \in \mathbb{R}$ which is contained in the Bohl dichotomy resolvent set of \eqref{1}
\begin{equation*}
   \rho_{\mathrm{BD}}(A)
   \coloneqq
   \mathbb{R} \setminus \Sigma _{\mathrm{BD}}(A)
   =
   \left\{ \gamma \in \mathbb{R}:\eqref{gamma-shift}%
    \text{ has a Bohl dichotomy}\right\}
\end{equation*}
or the exponential dichotomy resolvent set
\begin{equation*}
   \rho_{\mathrm{ED}}(A)
   \coloneqq
   \mathbb{R} \setminus \Sigma _{\mathrm{ED}}(A)
   =
   \left\{ \gamma \in \mathbb{R}:\eqref{gamma-shift}%
    \text{ has an exponential dichotomy}\right\}
\end{equation*}
the $\gamma$-shifted system \eqref{gamma-shift} admits a dichotomy on a splitting $(L_1,L_2)$. The subspace $L_1$ and, in particular, its dimension $k := \dim L_1$ is unique for a fixed $\gamma$ in the resolvent set (Proposition \ref{prop:uni_rep_dicho_spaces}).

Covers \(\mathcal U_1 \subseteq  \mathcal G(L_1)\), \(\mathcal U_2 \subseteq  \mathcal G(L_2)\) of $L_1$ and $L_2$, respectively, on which the dichotomy estimates hold uniformly (cf.\ Definition \ref{def:dicho}) are
\begin{equation*}
   \mathcal U_1
   \coloneqq
   \mathcal G_1(L_1)
   , \quad
   \mathcal U_2
   \coloneqq
   \mathcal G_1(L_2)
   \qquad \text{if }
   \gamma \in \rho_{\mathrm{BD}}(A)
\end{equation*}
and
\begin{equation*}
   \mathcal U_1
   \coloneqq
   \mathcal G_{k}(L_1)
   =
   \{L_1\}
   , \quad
   \mathcal U_2
   \coloneqq
   \mathcal G_{d-k}(L_2)
   =
   \{L_2\}
   \qquad \text{if }
   \gamma \in \rho_{\mathrm{ED}}(A)
   .
\end{equation*}
In both cases
\begin{equation}\label{uniformity-spaces}
   \mathcal U_1
   =
   \mathcal G_{j_1}(L_1)
   , \quad
   \mathcal U_2
   =
   \mathcal G_{j_2}(L_2)
\end{equation}
with $(j_1,j_2) \in \{(1,1), (k,d-k)\}$.
We extend this idea and answer in the following proposition for \(k \in \{0,\dots,d\}\) the question which \((j_1,j_2) \in \N\times\N\) have the property that for each splitting $(L_1,L_2)$ with $\dim L_1 = k$ the spaces $\mathcal U_1$ and $\mathcal U_2$ in \eqref{uniformity-spaces} are covers of $L_1$ and $L_2$, respectively (cf.\ Definition \ref{def:dicho}).

\begin{definition}[$k$-admissible pair]
\label{def:unif-dim}
Let $k, j_1, j_2 \in \{0, \dots, d\}$. We say that $(j_1, j_2)$ is \emph{$k$-admissible} if one of the following holds:

(a) If $k = 0$, then $j_1 = 0$ and $j_2 \in \{1, \dots, d\}$.

(b) If $k \in \{1, \dots, d -1\}$, then $j_1 \in \{1, \dots, k\}$ and $j_2 \in \{1, \dots, d-k\}$.

(c) If $k=d$, then $j_1 \in \{1, \dots, d\}$ and $j_2 = 0$
\end{definition}

\begin{proposition}[$k$-admissible pairs and covers of splittings]
\label{prop:unif-dim}
Let \(k, j_1, j_2 \in \{0,\dots,d\}\). Then $(j_1, j_2)$ is $k$-admissible if and only if for each splitting \((L_1, L_2)\) with $\dim L_1=k$, the spaces $L_1$ and $L_2$ are covered by their $j_1$- and $j_2$-dimensional subspaces, respectively, i.e.
\begin{equation*}
    \bigcup_{U\in\mathcal G_{j_1}(L_1)}U = L_1
    \qquad\text{and}\qquad
    \bigcup_{U\in\mathcal G_{j_2}(L_2)}U = L_2.
\end{equation*}
\end{proposition}

\begin{proof}
This follows from \(\mathcal G_j(L) = \emptyset\) for $L \in \mathcal G(\mathbb{R}^d)$ if \(\dim L < j\).
\end{proof}

\begin{definition}[Dichotomy with uniformity dimensions]
\label{def:dich-unif-dim}
Let \((j_1,j_2) \in \N\times\N\). We say that system \eqref{1} has a \emph{dichotomy with uniformity dimensions \((j_1,j_2)\)} on \((L_1, L_2)\) if system \eqref{1} has a dichotomy on \((L_1, L_2)\) uniformly on all subspaces of dimensions \((j_1,j_2)\), i.e.\ with
\begin{equation*}
   \mathcal U_1
   =
   \mathcal{G}_{j_1}(L_1)
   \qquad \text{and} \qquad
   \mathcal U_2
   =
   \mathcal{G}_{j_2}(L_2).
\end{equation*}
\end{definition}

The following remark is a consequence of the refinement Theorem \ref{thm:dicho-refine} applied to dichotomies which are uniform on all subspaces of prescribed dimensions \((j_1,j_2)\).

\begin{remark}[Refining uniformity dimensions]
\label{rem:dich-unif-dim}
If system \eqref{1} has a dichotomy with uniformity dimensions \((j_1,j_2)\) on \((L_1, L_2)\), then $(j_1,j_2)$ is $k$-admissible with $k = \dim L_1$ and system \eqref{1} has also a dichotomy with uniformity dimensions \((\ell_1,\ell_2)\) for each $k$-admissible \((\ell_1,\ell_2)\) which satisfies \(\ell_1 \leq j_1\) and \(\ell_2 \leq  j_2\).
\end{remark}

We are now ready to define a notion of spectrum based on dichotomies with $k$-admissible uniformity dimensions \((j_{1k},j_{2k})\) for the \(\gamma\)-shifted system \eqref{gamma-shift} where $k$ is the dimension of the dynamically characterized set (Proposition \ref{prop:uni_rep_dicho_spaces})
\begin{equation*}
    \{x_0 \in \R^d :
    \lim_{n\to\infty}\mathrm e^{-\gamma n}x(n,x_0) = 0
    \}.
\end{equation*}
The exponential and Bohl dichotomy spectrum are special cases for specific choices of uniformity dimensions (Remark \ref{rem:admissible_Bohl_exp}).

\begin{definition}[Dichotomy resolvent and spectrum]
\label{def:dichotomy-spectrum}
For each \(k \in \{0,\dots,d\}\) let \((j_{1k},j_{2k})\) be \(k\)-admissible. We call $J \coloneqq \big((j_{10},j_{20}), \dots, (j_{1d},j_{2d})\big) \in (\N\times\N)^{d+1}$ \emph{admissible uniformity dimensions}.
The \emph{dichotomy resolvent of \eqref{1} with admissible uniformity dimensions} $J$ is defined as
\begin{align*}
    \rho_J(A)
    \coloneqq
    &\{\gamma \in \R \,\vert\, \text{the \(\gamma\)-shifted system has a dichotomy on a splitting \((L_1,L_2)\),}
    \\
    &\quad\text{uniform on subspaces with dimension \((j_{1k},j_{2k})\), where \(k = \dim L_1\)}\}
\end{align*}
and the \emph{dichotomy spectrum of system \eqref{1} with admissible uniformity dimensions $J$ is}
\begin{equation*}
    \Sigma_J(A)
    \coloneqq
    \R \setminus \rho_J(A).
\end{equation*}
\end{definition}

\begin{remark}[Special cases of Bohl and exponential dichotomy spectrum]\label{rem:admissible_Bohl_exp}\hfill
Using Remark \ref{rem:uni_sub_Bohl_expo} on the uniformity subspaces of Bohl and exponential dichotomies it follows that

(a) \(\Sigma_\mathrm{BD}(A) = \Sigma_J(A)\) for
\(J = \big((0,1),(1,1),\dots,(1,1),(1,0)\big)\).

(b) \(\Sigma_\mathrm{ED}(A) = \Sigma_J(A)\) for
\(J = \big((0,d),(1,d-1),\dots,(d,0)\big)\).

In \cite[Remark 35]{Czornik2023} it was shown that in general $\Sigma_{\mathrm{BD}}(A) \subsetneq \Sigma_{\mathrm{ED}}(A)$. For an explicit construction of an $A$ with this strict inclusion, see e.g.\ \cite{Barreira2018b, Czornik2019}.
\end{remark}

\begin{remark}[Resolvent set characterized by dichotomy estimates]\label{rem:equiv_res_set}
For admissible uniformity dimensions $J$ the following two statements are equivalent:

(i) \(\gamma \in \rho_J(A)\),

(ii) there exists a splitting \((L_1,L_2)\) of \(\R^d\) and \(\alpha>0\), such that with \(k \coloneqq \dim L_1\) the dichotomy estimate \(D_1(\gamma-\alpha,U)\) holds uniformly for all \(U \in \mathcal G_{j_{1k}}(L_1)\) and \(D_2(\gamma+\alpha,U)\) holds uniformly for all \(U \in \mathcal G_{j_{2k}}(L_2)\).
\end{remark}

We now define exponents, called \emph{limiting Bohl exponents}, which turn out to be the boundary points of the dichotomy spectrum. Such descriptions were previously known
for the uniform exponential dichotomy only in the one-dimensional case \cite{Poetzsche2018} and in the multidimensional case only for the largest and
smallest elements of the spectrum \cite[Remark 23]{Czornik2023}.

\begin{definition}[Limiting Bohl exponents]
\label{def:limiting-Bohl}
For \(j,k\in\N\) define
\begin{align*}
    \overline\beta_{k,j} &\coloneqq \inf_{L\in\mathcal G_{k}(\R^d)}\sup_{U\in\mathcal G_j(L)}\overline\beta(U),
    \\
    \underline\beta_{k,j} &\coloneqq \sup_{L\in\mathcal G_{d-k}(\R^d)}\inf_{U\in\mathcal G_j(L)}\underline\beta(U).
\end{align*}
\end{definition}

\begin{remark}\label{rem:extremely_limiting}
Let \(k \in \{0,\dots,d\}\) and \((j_1,j_2)\) be \(k\)-admissible.
From Remark \ref{rem:infty_Bohl_exponents} it follows that for \(k = 0\), resp.\ \(k = d\),
\begin{equation*}
    \overline\beta_{k,j_1} = \overline\beta_{0,0} = -\infty
    \qquad\text{resp.}\qquad
    \underline\beta_{k,j_2} = \underline\beta_{d,0} = \infty.
\end{equation*}
From Remark \ref{rem:infty_Bohl_exponents} it also follows that
\begin{align*}
    &\overline\beta_{k,j_1} \in \big[-\ln\Vert A^{-1}\Vert_\infty,\ln\Vert A\Vert_\infty\big], && k \in \{1,\dots,d\},
    \\
    &\underline\beta_{k,j_2} \in \big[-\ln\Vert A^{-1}\Vert_\infty,\ln\Vert A\Vert_\infty\big], && k \in \{0,\dots,d-1\}.
\end{align*}
\end{remark}

\begin{lemma}[Characterization of dichotomy]\label{lem:res_int}
Let \(k \in \{0,\dots,d\}\) and \((j_1,j_2)\) be \(k\)-admissible. Then the following two statements are equivalent:

(i) \(\gamma \in \big(\overline\beta_{k,j_1},\underline\beta_{k,j_2}\big)\),

(ii) there exists a splitting \((L_1,L_2)\) of \(\R^d\) with \(\dim L_1 = k\), such that the \(\gamma\)-shifted system \eqref{gamma-shift} has a dichotomy on \((L_1,L_2)\) with uniformity dimensions \((j_1, j_2)\).

If (i) and (ii) hold, then $L_1$ is independent of  \(\gamma \in \big(\overline\beta_{k,j_1},\underline\beta_{k,j_2}\big)\). In case $\mathbb{T} = \mathbb{Z}$, also $L_2$ is independent of  \(\gamma \in \big(\overline\beta_{k,j_1},\underline\beta_{k,j_2}\big)\).
\end{lemma}

\begin{proof}
(i) $\Rightarrow$ (ii).
Let \(\gamma \in \big(\overline\beta_{k,j_1},\underline\beta_{k,j_2}\big)\).
There exists \(\alpha > 0\) with \(\gamma - \alpha > \overline\beta_{k,j_1}\) and \(\gamma + \alpha < \underline\beta_{k,j_2}\).
Applying Proposition \ref{prop:Rel_upper}\ref{prop:Rel_upper_2} (where \(\mathcal V\) is set to \(\mathcal G_k(\R^d)\) and \(\mathcal U(V)\) is set to \(\mathcal G_{j_1}(V)\) for \(V \in \mathcal V\)), we obtain the existence of a space \(L_1 \in \mathcal G_k(\R^d)\), such that \(D_1(\gamma-\alpha,U)\) holds uniformly for all \(U \in \mathcal G_{j_1}(L_1)\).
Applying Proposition \ref{prop:Rel_lower}\ref{prop:Rel_lower_2} (where \(\mathcal V\) is set to \(\mathcal G_{d-k}(\R^d)\) and \(\mathcal U(V)\) is set to \(\mathcal G_{j_2}(V)\) for \(V \in \mathcal V\)), we obtain the existence \(L_2 \in \mathcal G_{d-k}(\R^d)\), such that \(D_2(\gamma+\alpha,U)\) holds uniformly for all \(U \in \mathcal G_{j_2}(L_2)\).
From Proposition \ref{prop:est_dicho}, we conclude \(L_1 \cap L_2 = \{0\}\), so that \((L_1,L_2)\) is a splitting of \(\R^d\) with \(\dim L_1 = k\). With Remark \ref{rem:equiv_res_set} (ii) follows.

(ii) $\Rightarrow$ (i).
Suppose that \((L_1,L_2)\) is a splitting of \(\R^d\) with \(\dim L_1 = k\), such that the \(\gamma\)-shifted system \eqref{gamma-shift} has a dichotomy on \((L_1,L_2)\) with uniformity dimensions \((j_1, j_2)\).
By Remark \ref{rem:equiv_res_set}, there is \(\alpha > 0\), such that \(D_1(\gamma-\alpha,U)\) holds uniformly for all \(U \in \mathcal G_{j_1}(L_1)\) and \(D_2(\gamma+\alpha,U)\) holds uniformly for all \(U \in \mathcal G_{j_2}(L_2)\) for some \(\alpha > 0\).
Propositions \ref{prop:Rel_upper}\ref{prop:Rel_upper_1} and \ref{prop:Rel_lower}\ref{prop:Rel_lower_1} prove that \(\gamma \in \big(\overline\beta_{k,j_{1k}},\underline\beta_{k,j_{2k}}\big)\).

To show that $L_1$ is independent of  \(\gamma \in \big(\overline\beta_{k,j_1},\underline\beta_{k,j_2}\big)\), let  \(\gamma, \gamma' \in \big(\overline\beta_{k,j_1},\underline\beta_{k,j_2}\big)\) with $\gamma < \gamma'$ and let \((L_1,L_2)\) and \((L_1',L_2')\) be the respective splittings. By Proposition \ref{prop:uni_rep_dicho_spaces}
\begin{equation*}
    L_1 = \{x_0 \in \R^d : \lim_{n\to\infty}\mathrm e^{-\gamma n}x(n,x_0) = 0\},
\end{equation*}
and
\begin{equation*}
    L_1' = \{x_0 \in \R^d : \lim_{n\to\infty}\mathrm e^{-\gamma' n}x(n,x_0) = 0\}.
\end{equation*}
Since $\gamma < \gamma'$, it follows that $L_1 \subseteq L_1'$. Using the fact that \(\dim L_1 = \dim L_1' = k\), we conclude that $L_1 = L_1'$. Similarly it follows that $L_2 = L_2'$ in case $\mathbb{T} = \mathbb{Z}$.
\end{proof}

\begin{lemma}[Dichotomy resolvent]\label{lem:dicho_res}
Let $J \coloneqq \big((j_{10},j_{20}), \dots, (j_{1d},j_{2d})\big) \in (\N\times\N)^{d+1}$ be admissible uniformity dimensions (in the sense of Definition \ref{def:dichotomy-spectrum}). Then
\begin{align*}
    \rho_J(A)
    &= \big(-\infty, \underline\beta_{0,j_{20}}\big) \cup {}
\\
    &\phantom{i} \quad
    \big(\overline\beta_{1,j_{11}},
    \underline\beta_{1,j_{21}}\big)
    \cup \dots \cup
    \big(\overline\beta_{(d-1),j_{1(d-1)}},
    \underline\beta_{(d-1),j_{2(d-1)}}\big)
    \cup {}
\\
    &\phantom{i} \quad
    \big(\overline\beta_{d,j_{1d}}, \infty\big)
\end{align*}
and the union is ordered with
\begin{equation*}
    -\ln\Vert A^{-1}\Vert_\infty \leq \underline\beta_{k,j_{2k}} \leq \overline\beta_{k+1,j_{1(k+1)}} \leq \ln\Vert A\Vert_\infty,
    \qquad k \in \{0,\dots, d-1\}.
\end{equation*}
\end{lemma}

\begin{proof}
From Lemma \ref{lem:res_int} we obtain
\begin{equation*}
    \rho_J(A) = \bigcup_{k=0}^d\big(\overline\beta_{k,j_{1k}}, \underline\beta_{k,j_{2k}}\big).
\end{equation*}
The equality follows from Remark \ref{rem:extremely_limiting}.
We show for \(k \in \{0,\dots, d-1\}\) that \(\underline\beta_{k,j_{2k}} \leq \overline\beta_{k+1,j_{1(k+1)}}\), i.e.
\begin{equation*}
    \sup_{L\in\mathcal G_{d-k}(\R^d)}\inf_{U\in\mathcal G_{j_{2k}}(L)}\underline\beta(U)
    \leq
    \inf_{L\in\mathcal G_{k+1}(\R^d)}\sup_{U\in\mathcal G_{j_{1(k+1)}}(L)}\overline\beta(U).
\end{equation*}
To this end, let \(L \in \mathcal G_{d-k}(\R^d)\), \(\widetilde L \in \mathcal G_{k+1}(\R^d)\).
Since \(\dim L = d-k\), \(\dim\widetilde L = k+1\), there is \(x_0 \in L \cap \widetilde L\) with \(x_0 \neq 0\).
Since \(j_{2k}, j_{1(k+1)} > 0\) by admissibility, there is \(V \in \mathcal G_{j_{2k}}(L)\), \(\widetilde V \in \mathcal G_{j_{1(k+1)}}(\widetilde L)\) with \(x_0 \in V \cap \widetilde V\).
Thus by Remark \ref{rem:infty_Bohl_exponents},
\begin{equation*}
    \inf_{U\in\mathcal G_{j_{2k}}(L)}\underline\beta(U)
    \leq
    \underline\beta(V)
    \leq
    \underline\beta(\operatorname{span}\{x_0\})
    \leq
    \overline\beta(\operatorname{span}\{x_0\})
    \leq
    \underline\beta(\widetilde V)
    \leq
    \sup_{U\in\mathcal G_{j_{1(k+1)}}(\widetilde L)}\underline\beta(U).
\end{equation*}
\end{proof}

We are now in a position to formulate our first main result on the dichotomy spectrum which generalizes Theorem \ref{thm:exp-bohl-spectrum} and provides a spectral flag or spectral filtration of dynamically characterized subspaces which are also related to the formulas for the spectrum given in our second main result Theorem \ref{thm:spectral_formula} below.

\begin{theorem}[Spectral theorem]\label{thm:spectral}
Let \(J \in (\N\times\N)^{d+1}\) be admissible.
The dichotomy spectrum \(\Sigma_J(A)\) is the nonempty union of at most \(d\) compact intervals
\begin{equation*}
    \Sigma_J(A)
    =
    [a_1,b_1] \cup \dots \cup [a_\ell,b_\ell],
\end{equation*}
where \(a_1 \leq b_1 < a_2 \leq b_2 < \dots < a_\ell \leq b_\ell\) and \(\ell \in \{1,\dots,d\}\).
Setting \(b_0 \coloneqq -\infty\) and \(a_{\ell+1} \coloneqq \infty\), there exists a spectral filtration
\begin{equation*}
    \{0\} = M_0 \subsetneq M_1 \subsetneq \dots \subsetneq M_\ell = \R^d,
\end{equation*}
of subspaces \(M_k \subseteq \R^d\), \(k \in \{0,\dots,\ell\}\) with
\begin{equation*}
    M_k = \big\{x_0 \in \R^d \,\vert\, \lim_{n\to\infty}\mathrm e^{-\gamma n}x(n,x_0) = 0\big\},
    \qquad
    \gamma \in (b_k,a_{k+1}).
\end{equation*}
For two-sided time \(\T = \Z\), there exists a spectral decomposition
\begin{equation*}
    \R^d = W_1 \oplus \dots \oplus W_\ell
\end{equation*}
into subspaces \(W_k \subseteq \R^d\), with
\begin{equation*}
    W_k
    =
    \big\{x_0 \in \R^d \,\vert\, \lim_{n\to-\infty}\mathrm e^{-\gamma_1n}x(n,x_0) = \lim_{n\to\infty}\mathrm e^{-\gamma_2n}x(n,x_0) = 0\big\},
\end{equation*}
where \(\gamma_1 \in (b_{k-1},a_k)\), \(\gamma_2 \in (b_k,a_{k+1})\) and \(k \in \{1,\dots,\ell\}\).
\end{theorem}

\begin{proof}
We define \(k_0 \coloneqq 0\) and iteratively if \(k_i\) is defined, we define \(k_{i+1} \coloneqq \inf\{n \in \{k_i+1,\dots,d\} : (\overline\beta_{n,j_{1n}},\underline\beta_{n,j_{2n}}) \neq \emptyset\}\), \(i \in \N\).
We define \(\ell \in \{1,\dots,d\}\) as the maximal number, such that \(k_\ell \neq \infty\).
Then \(k_\ell = d\) and we set
\begin{align*}
    a_i \coloneqq \underline\beta_{k_{i-1},j_{2k_{i-1}}}, \qquad i \in \{1,\dots,\ell+1\}
    \quad\text{and}\quad
    b_i \coloneqq \overline\beta_{k_i,j_{1k_i}}, && i \in \{0,\dots,\ell\}.
\end{align*}
By Lemma \ref{lem:dicho_res}, we obtain \(a_1 \leq b_1 < a_2 \leq b_2 < \dots < a_\ell \leq b_\ell\) and
\begin{equation*}
    \Sigma_J(A)
    =
    [a_1,b_1] \cup \dots \cup [a_\ell,b_\ell].
\end{equation*}
To show the existence of the spectral filtration, choose for each \(k \in \{0,\dots,\ell\}\) a $\gamma \in (b_k,a_{k+1})$ and define
\[
    M_k 
    \coloneqq
    \big\{x_0 \in \R^d \,\vert\, \lim_{n\to\infty}\mathrm e^{-\gamma n}x(n,x_0) = 0\big\}.
\]
By Proposition \ref{prop:uni_rep_dicho_spaces} the $\gamma$-shifted system has a dichotomy on $(L_1,L_2)$ with $L_1 = M_k$. The fact that $M_k$ ist well-defined, i.e.\ independent of $\gamma \in (b_k,a_{k+1})$, follows from Proposition \ref{lem:res_int}. Obviously $M_k \subseteq M_{k+1}$ for \(k \in \{0,\dots,\ell-1\}\). By Proposition \ref{lem:res_int} $\dim M_k < \dim M_{k+1}$, proving that the inclusion is strict.

Now suppose \(\T = \Z\).
Then by Proposition \ref{prop:uni_rep_dicho_spaces} and Lemma \ref{lem:res_int} the space
\begin{equation*}
    N_k \coloneqq \big\{x_0 \in \R^d : \lim_{n\to-\infty}\mathrm e^{-\gamma n}x(n,x_0) = 0\big\},
    \qquad
    \gamma \in (b_k,a_{k+1}),
\end{equation*}
is well-defined for \(k \in \{0,\dots,\ell\}\), so that \(W_k = M_k \cap N_{k-1}\) for \(k \in \{1,\dots,\ell\}\). By Proposition \ref{prop:uni_rep_dicho_spaces} the $\gamma$-shifted system has a dichotomy on $(L_1,L_2)$ with $L_1 = M_k$ and $L_2 = N_k$, implying that $M_k \cap N_k = \{0\}$ and $M_k + N_k = \mathbb{R}^d$.

We now show the spectral decomposition.
To show that the sum is direct, let \(\ell,j\in\{1,\dots,d\}\) with \(\ell < j\) and let \(x_0 \in (M_\ell \cap N_{\ell-1}) \cap (M_j \cap N_{j-1})\).
Since \(M_\ell \subseteq M_{j-1}\), we have \(x_0 \in M_{j-1} \cap N_{j-1} = \{0\}\) by Proposition \ref{prop:est_dicho}, proving that \(x_0 = 0\).
We show inductively that
\begin{align*}
    &(a) && M_1 = M_1\cap N_0.
    \\
    &(b) && \text{If \(k > 1\) and \(M_{k-1} = \bigoplus_{j = 1}^{k-1} M_j\cap N_{j-1}\) for \(k \in \{1,\dots,d\}\)},
    \\
    & && \qquad\text{then \(M_k = \bigoplus_{j = 1}^k M_j\cap N_{j-1}\)}.
\end{align*}
Part (a) follows, since \(N_0 = \R^d\).
To show part (b), suppose that \(k > 1\) and \(M_{k-1} = \bigoplus_{j = 1}^{k-1} M_j\cap N_{j-1}\).
That \(\bigoplus_{j = 1}^k M_j\cap N_{j-1} \subseteq M_k\), follows from \(M_j \subseteq M_k\) for \(j \in \{1,\dots,k\}\).
Now suppose that \(x_0 \in M_k\).
Let \(x_0 = x_1 + x_2\), with \(x_1 \in M_{k-1}\) and \(x_2 \in N_{k-1}\).
Since \(x_0 - x_1 \in M_k\), we have \(x_2 = x_0 - x_1 \in M_k \cap N_{k-1}\).
By assumption also \(x_1 \in M_{k-1} = \bigoplus_{j = 1}^{k-1} M_j\cap N_{j-1}\), so that \(x_0 = x_1 + x_2 \in \bigoplus_{j = 1}^k M_j\cap N_{j-1}\), indeed.
The theorem now follows, since \(M_\ell = \R^d\).
\end{proof}

The proof of Theorem \ref{thm:spectral} yields explicit formulas for the boundary points of the spectral intervals.

\begin{theorem}[Formula for the dichotomy spectrum]\label{thm:spectral_formula}
Let \(J \in (\N\times\N)^{d+1}\) be admissible.
Let \(k_i \coloneqq \dim M_i\) denote the dimension of the spectral flag subspaces \(M_i\), \(i \in \{0,\dots,\ell\}\), of Theorem \ref{thm:spectral}.
The following formula holds for the dichotomy spectrum \(\Sigma_J(A) = [a_1,b_1] \cup \dots \cup [a_\ell,b_\ell]\) of system \eqref{1}.
\begin{equation*}
    a_i = \underline\beta_{k_{i-1},j_{2k_{i-1}}}
    \qquad\text{and}\qquad
    b_i = \overline\beta_{k_i,j_{1k_i}},
    \qquad
    i \in \{1,\dots,\ell\}.
\end{equation*}
\end{theorem}

\begin{remark}[Formulas for the Bohl and exponential dichotomy spectrum]
Using the notation of Theorem \ref{thm:spectral_formula}, we obtain the following formulas for the boundary points of the Bohl dichotomy spectrum and the exponential dichotomy spectrum (cf.\ Remark \ref{rem:admissible_Bohl_exp}).

(a) The boundary points of the Bohl dichotomy spectrum are
\begin{equation*}
    a_i = \underline\beta_{k_{i-1},1}
    \qquad\text{and}\qquad
    b_i = \overline\beta_{k_i,1},
    \qquad
    i \in \{1,\dots,\ell\}.
\end{equation*}

(b) The boundary points of the exponential dichotomy spectrum are
\begin{equation*}
    a_i = \underline\beta_{k_{i-1},d-k_{i-1}}
    \qquad\text{and}\qquad
    b_i = \overline\beta_{k_i,d-k_i},
    \qquad
    i \in \{1,\dots,\ell\}.
\end{equation*}
\end{remark}

In case system \eqref{1} has a dichotomy then the limiting Bohl exponents of Definition \ref{def:limiting-Bohl} can also be expressed in terms of the splitting of the dichotomy.

\begin{corollary}[Limiting Bohl exponents in case of dichotomy]
If system \eqref{1} has a dichotomy on \((L_1,L_2)\) with uniformity dimensions \((j_1, j_2)\) then
\begin{equation}\label{col:of_main:eq:formula_1}
    \overline\beta_{k,j_1} = \sup_{U\in\mathcal G_{j_1}(L_1)}\overline\beta(U)
\end{equation}
where $k = \dim L_1$,
and if \(\T = \Z\), then
\begin{equation}\label{col:of_main:eq:formula_2}
    \underline\beta_{k,j_2} = \inf_{U\in\mathcal G_{j_2}(L_2)}\overline\beta(U).
\end{equation}
\end{corollary}

\begin{proof}
The first part follows readily from Lemma \ref{lem:res_int} and we only proof formula \eqref{col:of_main:eq:formula_1}.
Suppose that system \eqref{1} has a dichotomy with splitting \((L_1,L_2)\) uniform on subspaces of dimensions \(j_1\) and \(j_2\).
Let \(\alpha>0\) with \(-\alpha \in (\overline\beta_{k,j_1},0)\), so that \(\mathrm D_1(-\alpha,U)\) holds uniformly for all \(U \in \mathcal G_{j_1}(L_1)\) by Proposition \ref{prop:Rel_upper}\ref{prop:Rel_upper_2}.
Let \(k \coloneqq \dim L_{1}\).
If \(k = 0\), then \(j_1 = 0\) and \eqref{col:of_main:eq:formula_1} is satisfied (cf.\ Remark \ref{rem:infty_Bohl_exponents}).
Let \(k > 0\).
Suppose that \(L_{1}' \in \mathcal{G}_{k}(\mathbb{R}^{d})\) with
\begin{equation}\label{col:of_main:eq:formula_1_proof}
    \sup_{U' \in \mathcal{G}_{j_{1}}(L_{1}')}\overline\beta(U')
    \leq
    \sup_{U \in \mathcal{G}_{j_{1}}(L_{1})}\overline\beta(U).
\end{equation}
It holds that \(\overline\beta(U) \leq -\alpha\) for all \(U \in \mathcal G_{j_1}(L_1)\) (cf.\ Proposition \ref{prop:Rel_upper}\ref{prop:Rel_upper_1} with \(\mathcal V \coloneqq \{L_1\}\) and \(\mathcal U(L_1) \coloneqq \mathcal G_{j_1}(L_1)\)).
By the inequality \eqref{col:of_main:eq:formula_1_proof}, also \(\overline\beta(U') \leq -\alpha\) for all \(U' \in \mathcal G_{j_1}(L_1')\).
By Proposition \ref{prop:Rel_upper}\ref{prop:Rel_upper_2}, we conclude that \(\mathrm D_1(-\alpha/2,U')\) holds for all \(U' \in \mathcal G_{j_1}(L_1')\), so that
\begin{equation*}
    \underset{n\rightarrow \infty }{\lim }x(n,x_{0})=0,
    \qquad
    x_{0} \in L_{1}'.
\end{equation*}
This implies that \(L_{1}' \subseteq L_{1}\) by Proposition \ref{prop:uni_rep_dicho_spaces}.
But since \(\dim L_{1} = \dim L_{1}'\), also \(L_{1}=L_{1}'\).
We conclude that for all \(L_{1}' \in \mathcal{G}_{k}(\mathbb{R}^{d})\), we have
\begin{equation*}
    \sup_{U' \in \mathcal{G}_{j_{1}}(L_{1}')}\overline\beta(U')
    \geq
    \sup_{U \in \mathcal{G}_{j_{1}}(L_{1})}\overline\beta(U).
\end{equation*}
Taking the infimum over all \(L_{1}' \in \mathcal{G}_{k}(\mathbb{R}^{d})\) yields formula \eqref{col:of_main:eq:formula_1}.
\end{proof}

\begin{remark}[Open problem]
It is not clear to the authors, if \eqref{col:of_main:eq:formula_2} holds for all \(L_2\), complementary to \(L_1\) in the case \(\T = \N\).
\end{remark}

\section{Maximal subspaces of uniformity}

In this section we introduce maximality of the uniformity dimensions of a dichotomous system in Definition \ref{def:max-unif-dimensions} and discuss its dependence on the dichotomy splitting in Theorems \ref{thm:uni_indep_L_2} and \ref{thm:max-unif-dim}. An open problem and a conjecture are formulated in Remark \ref{rem:max-unif-dim}.

\begin{definition}[Maximal uniformity dimensions]
\label{def:max-unif-dimensions}
Suppose system \eqref{1} has a dichotomy on \((L_1, L_2)\). Then $(u_1,u_2)$ are called \emph{maximal uniformity dimensions} if

(a) system \eqref{1} has a dichotomy on \((L_1, L_2)\) with uniformity dimensions \((u_1,u_2)\), and

(b) if \eqref{1} has a dichotomy on \((L_1, L_2)\) with uniformity dimensions \((\ell_1, \ell_2)\) then
\begin{equation*}
    \ell_1 \in \{0,\dots,u_1\}
    \qquad\text{and}\qquad
    \ell_2 \in \{0,\dots,u_2\}.
\end{equation*}
\end{definition}

\begin{remark}[Dichotomy implies existence of maximal uniformity dimensions]
\label{rem:max-unif-dimensions}
If system \eqref{1} has a dichotomy on \((L_1, L_2)\) there exist maximal uniformity dimensions \((u_1, u_2) \in \{0,\dots,\dim L_1\} \times \{0,\dots,\dim L_2\}\).
\end{remark}

In case system \eqref{1} has a dichotomy on \((L_1, L_2)\) then by Proposition \ref{prop:uni_rep_dicho_spaces} the space \(L_1\) is unique and hence the maximal uniformity dimension \(u_1\) does not depend on \(L_1\). If \(\T = \Z\), the space \(L_2\) is unique and in this case the maximal uniformity dimension \(u_2\) does not depend on \(L_2\). We provide a partial answer to the question how \(u_2\) does depend on \(L_2\) in Theorem \ref{thm:max-unif-dim} and formulate an open problem in Remark \ref{rem:max-unif-dim}. As preparatory results we prove the following Lemma and Theorem \ref{thm:uni_indep_L_2}.

\begin{lemma}[Dichotomy estimate follows from uniform estimate at later time]\label{lem:Adams_L_2_thoughts}
Let $\gamma \in \R$, \(U \in \mathcal G(\R^d)\), \(C>0\) and \(m_0 \in \N\) with
\begin{equation}\label{lem:Adams_L_2_thoughts_1}
    \Vert x(n,x_0)\Vert \geq C \mathrm e^{\gamma(n-m)}\Vert x(m,x_0)\Vert,
    \qquad m,n\in\N,\,n\geq m\geq m_0,\, x_0 \in U.
\end{equation}
Then \(\mathrm D_2(\gamma,U)\) holds uniformly.
\end{lemma}

\begin{proof}
If \(m_0 = 0\), then the statement of Lemma \ref{lem:Adams_L_2_thoughts} is clear.
Suppose that \(m_0 \geq 1\).
By induction, it suffices to show that
\begin{equation}\label{lem:Adams_L_2_thoughts_2}
    \Vert x(n,x_0)\Vert \geq C'\mathrm e^{\gamma (n-m_0+1)}\Vert
    x(m_0-1,x_0)\Vert,
    \qquad n\in \N,\, n\geq m_0-1,\,x_0 \in U
\end{equation}
holds for some \(C' > 0\).
Let us denote
\begin{equation*}
    a = \max\big\{\Vert A\Vert _\infty,\Vert
A^{-1}\Vert_\infty\big\}.
\end{equation*}
Considering \eqref{lem:Adams_L_2_thoughts_1} for \(n\) set to \(n+1\) and for \(m\) set to \(m_0\), we get
\begin{equation*}
    \Vert x(n+1,x_0)\Vert \geq C\mathrm e^{\gamma (n-m_{0}+1)}\Vert x(m_0,x_0)\Vert,
    \qquad m,n\in \N,\,n\geq m_0-1,\,x_0\in U
\end{equation*}
and consequently
\begin{equation*}
    \Vert A^{-1}(m_0-1)\Vert \cdot \Vert x(n+1,x_0)\Vert \geq C\mathrm e^{\gamma (n-m_0+1)}\Vert A^{-1}(m_0-1)\Vert \cdot \Vert x(m_0,x_0)\Vert.
\end{equation*}
Noting that
\begin{align*}
    \Vert A^{-1}(m_0-1)\Vert \Vert x(n+1,x_0)\Vert  &\leq a\Vert
    x(n+1,x_{0})\Vert
    \\
    &= a\Vert A(n)x(n,x_0)\Vert
    \\
    &\leq a^{2}\Vert x(n,x_0)\Vert
\end{align*}
and
\begin{equation*}
    \Vert A^{-1}(m_0-1)\Vert \Vert x(m_{0},x_0)\Vert \geq \Vert A^{-1}(m_0-1)x(m_0,x_0)\Vert =\Vert x(m_0-1,x_0)\Vert,
\end{equation*}
we obtain
\begin{equation*}
    a^2\Vert x(n,x_0)\Vert \geq C\mathrm e^{\gamma(n-m_0+1)}\Vert x(m_0-1,x_0)\Vert .
\end{equation*}
This proves inequality \eqref{lem:Adams_L_2_thoughts_1} with \(C' = C/a^2\).
\end{proof}

\begin{theorem}[Independence of maximal uniformity dimensions of splitting]
\label{thm:uni_indep_L_2}
Suppose that \(T = \N\). Let \((L_{1},L_{2})\) and \((L_{1},L_{2}')\) be splittings of \(\R^d\).
Suppose that System \eqref{1} has a dichotomy on \((L_{1},L_{2})\) with maximal uniformity dimensions \((u_{1},u_{2})\) and on \((L_{1},L_{2}')\) and with maximal uniformity dimensions \((u_{1}',u_{2}')\).
Then
\begin{equation}\label{pr14_1}
    u_{1}=u_{1}'.
\end{equation}
If additionally for each \(V' \in \mathcal{G}_{u_{2}}(L_{2}')\), we have \(\dim V_{L_{1}}'\leq u_{1}\), then
\begin{equation}\label{pr14_2}
    u_{2}\leq u_{2}'.
\end{equation}
\end{theorem}

\begin{proof}
The equality \eqref{pr14_1} is clear.
Let \(\alpha>0\) be such that \(\mathrm D_{1}(-\alpha ,U)\) holds uniformly for all \(U\in \mathcal{G}_{j_{1}}(L_{1})\)
and \(\mathrm D_{2}(\alpha ,U)\) holds uniformly for all \(U\in \mathcal{G}_{j_{2}}(L_{2})\) for all \(j_{1}\in \{0,...,u_{1}\}\) and \(j_{2} \in \{0,...,u_{2}\}\).
Let \(V'\in \mathcal{G}_{u_{2}}(L_{2}')\).
To prove the inequality \eqref{pr14_2}, it suffices to show that the estimate \(\mathrm D_{2}(\alpha ,V')\) holds uniformly.
Since \(\dim \pi_{L_1}[V']\leq u_{1}\) and \(\dim \pi_{L_2}[V']\leq \dim V' = u_{2}\), there exist constants \(C_{1},C_{2}>0\) such that
\begin{equation}\label{pr14_3}
    \Vert x(n,x_{0})\Vert \leq C_{1}\mathrm{e}^{-\alpha (n-m)}\Vert x(m,x_{0})\Vert,
    \qquad m,n\in \N,\,n\geq m,\,x_{0}\in \pi_{L_1}[V']
\end{equation}
and
\begin{equation}\label{pr14_4}
    \Vert x(n,x_{0})\Vert \geq C_{2}\mathrm{e}^{\alpha (n-m)}\Vert x(m,x_{0})\Vert ,
    \qquad m,n\in \N,\,n\geq m,\,x_{0}\in \pi_{L_2}[V'].
\end{equation}
Let \(K > 0\) be such that
\begin{equation*}
    C_{1}\leq K \qquad\text{and}\qquad C_{2}\geq K^{-1}.
\end{equation*}
Then \eqref{pr14_3} and \eqref{pr14_4} imply that
\begin{equation}\label{pr14_5}
    \Vert x(n,x_{0})\Vert \leq K\mathrm{e}^{-\alpha (n-m)}\Vert x(m,x_{0})\Vert,
    \qquad m,n\in \N,\,n\geq m,\,x_{0}\in \pi_{L_1}[V']
\end{equation}
and
\begin{equation}\label{pr14_6}
    \Vert x(n,x_{0})\Vert \geq K^{-1}\mathrm{e}^{\alpha (n-m)}\Vert x(m,x_{0})\Vert,
    \qquad m,n\in \N,\,n\geq m,\,x_{0}\in \pi_{L_2}[V'].
\end{equation}
Consider \(v' \in V'\) with \(\Vert v'\Vert = 1\) and note that \(v' \neq 0\) implies \(\pi_{L_2}v' \neq 0\).
We have
\begin{equation*}
    \Vert x(n,v')\Vert \geq \Vert x(n,\pi_{L_2}v')\Vert - \Vert x(n,\pi_{L_1}v')\Vert , \qquad n \in \N
\end{equation*}
and by applying \eqref{pr14_5} and \eqref{pr14_6}, we get
\begin{align*}
    \Vert x(n,v')\Vert &\geq K^{-1}\mathrm{e}^{\alpha (n-m)}\Vert
    x(m,\pi_{L_2}v')\Vert -K\mathrm{e}^{-\alpha (n-m)}\Vert x(m,\pi_{L_1}v')\Vert
    \\
    &\geq K^{-1}\mathrm{e}^{\alpha (n-m)}\Vert x(m,\pi_{L_2}v')\Vert -K\mathrm{e}^{\alpha (n-m)}\Vert x(m,\pi_{L_1}v')\Vert
    \\
    &= \mathrm{e}^{\alpha (n-m)}( K^{-1}\Vert x(m,\pi_{L_2}v')\Vert -K\Vert
    x(m,\pi_{L_1}v')\Vert),
    \\
    &\mathrm{e}^{\alpha (n-m)}( K^{-1}-K\frac{\Vert x(m,\pi_{L_1}v')\Vert
    }{\Vert x(m,\pi_{L_2}v')\Vert}) \Vert x(m,\pi_{L_2}v')\Vert,
    \qquad m,n\in \N,\,n\geq m.
\end{align*}
Taking in \eqref{pr14_3} \(x_{0}=\pi_{L_1}v'\) and \(m=0\), and in \eqref{pr14_4} \(x_{0}=\pi_{L_2}v'\) and \(m = 0\), we get
\begin{equation}\label{pr14_6_5}
    \frac{\Vert x(k,\pi_{L_1}v')\Vert}{\Vert x(k,\pi_{L_2}v')\Vert}
    \leq
    \frac{K^{2}}{\mathrm{e}^{2\alpha k}}\frac{\Vert \pi_{L_1}v'\Vert}{\Vert \pi_{L_2}v'\Vert},
    \qquad k\in \N
\end{equation}
and therefore
\begin{equation}\label{pr14_7}
    \Vert x(n,v')\Vert \geq \mathrm{e}^{\alpha (n-m)}\bigg(K^{-1}-\frac{K^{3}}{\mathrm{e}^{2\alpha m}}\frac{\Vert \pi_{L_1}v'\Vert}{\Vert \pi_{L_2}v'\Vert}\bigg) \Vert x(m,\pi_{L_2}v')\Vert,
    \qquad m,n\in \N,\,n\geq m.
\end{equation}
The function \(v' \mapsto \frac{\Vert \pi_{L_1}v'\Vert}{\Vert \pi_{L_2}v'\Vert }\) defined on the compact set \(\{v' \in V : \Vert v'\Vert = 1\}\) is continuous and therefore attains its maximum value.
Hence
\begin{equation*}
    C' \coloneqq \max \bigg\{\frac{\Vert \pi_{L_1}v'\Vert}{\Vert \pi_{L_2}v'\Vert } : v'\in V,\,\Vert v'\Vert =1\bigg\} \in [0,\infty),
\end{equation*}
is well-defined.
The inequality \eqref{pr14_7} then implies that
\begin{equation*}
    \Vert x(n,v')\Vert \geq \mathrm{e}^{\alpha (n-m)}(K^{-1}-\frac{K^{3}}{\mathrm{e}^{2\alpha m}}C') \Vert x(m,\pi_{L_2}v')\Vert,
    \qquad m,n\in \N,\,n\geq m.
\end{equation*}
Consider \(m_{0} \in \N\) such that
\begin{equation*}
    K^{-1}-\frac{K^{3}}{\mathrm{e}^{2\alpha m_{0}}}C' >0.
\end{equation*}
Define
\begin{equation*}
    C'' \coloneqq K^{-1}-\frac{K^{3}}{\mathrm{e}^{2\alpha m_{0}}}C'.
\end{equation*}
We get
\begin{equation}\label{pr14_8}
    \Vert x(n,v')\Vert \geq \mathrm{e}^{\alpha (n-m)}C''\Vert x(m,\pi_{L_2}v')\Vert,
    \qquad m,n\in \N,\,n\geq m\geq m_{0},
\end{equation}
since
\begin{equation*}
    K^{-1}-\frac{K^{3}}{\mathrm{e}^{2\alpha m}}C'\geq C''
\end{equation*}
for \(m \geq m_{0}\).
Note that by \eqref{pr14_6_5}, we have
\begin{align*}
    \Vert x(k,v')\Vert  &\leq \Vert x(k,\pi_{L_2}v')\Vert +\Vert x(k,\pi_{L_1}v')\Vert
    \\
    &\leq \Vert x(k,\pi_{L_2}v')\Vert +\frac{K^{2}}{\mathrm{e}^{2\alpha k}}C'\Vert x(k,\pi_{L_2}v')\Vert
    \\
    &= ( 1+\frac{K^{2}}{\mathrm{e}^{2\alpha k}}C') \Vert x(k,\pi_{L_2}v')\Vert
    \\
    &\leq ( 1+K^{2}C') \Vert x(k,\pi_{L_2}v')\Vert,
    \qquad k\in \N
\end{align*}
and therefore
\begin{equation}\label{pr14_9}
    \Vert x(m,y(v'))\Vert \geq (1 + K^{2}C')^{-1}\Vert x(m,v')\Vert,
    \qquad m\in \N.
\end{equation}
From the inequalities \eqref{pr14_8} and \eqref{pr14_9}, we get
\begin{equation*}
    \Vert x(n,v')\Vert \geq \mathrm{e}^{\alpha (n-m)}C''( 1+K^{2}C') ^{-1}\Vert x(m,v')\Vert,
    \qquad m,n \in \N,\,n\geq m\geq m_{0}.
\end{equation*}
The statement now follows from Lemma \ref{lem:Adams_L_2_thoughts}.
\end{proof}

\begin{theorem}[Maximal uniformity dimensions for one-sided time]
\label{thm:max-unif-dim}
Suppose that \(T = \N\).
Let \((L_{1},L_{2})\) be a splitting of \(\R^d\).
Suppose that System \eqref{1} has a dichotomy on \((L_{1},L_{2})\) with maximal uniformity dimensions \((u_{1},u_{2})\).
Denote by \((u_1,u_2^\perp)\) the maximal uniformity dimensions of system \eqref{1} with respect to the splitting \((L_1,L_1^\perp)\).
Then it holds that \(u_2 \leq u_2^\perp\).
\end{theorem}

\begin{proof}
This follows from Theorem \ref{thm:uni_indep_L_2}, noting that for each \(V' \in \mathcal{G}_{u_{2}}(L_1^\perp)\), we have \(\dim V_{L_{1}}[V'] = 0 \leq u_{1}\)
\end{proof}

In the following remark we list those cases in which Theorem \ref{thm:uni_indep_L_2} yields a unique maximal uniformity dimension of a dichotomy, i.e.\ independent of the splitting, and we formulate a conjecture on the general dependence of the maximal uniformity dimensions on the subspace \(L_2\) which is complementary to \(L_1\).

\begin{remark}[Dependence of maximal uniformity dimensions on splitting for one-sided time]
\label{rem:max-unif-dim}
Suppose that system \eqref{1} with \(\T = \N\) has a dichotomy on a splitting $(L_1, L_2)$.

(a) The maximal uniformity dimensions do not depend on $L_2$ in case \(\dim L_1 = 1\) or \(\dim L_1 = d-1\).

(b) Conjecture: There exists \(A \colon \N \to \mathrm{GL}(\R^d)\) such that the maximal uniformity dimension \(u_2\) depends on the choice of $L_2$.
\end{remark}

\section*{Acknowledgements}

The research of Adam Czornik was supported by the Polish National Agency for Academic Exchange according to the decision No. BPN/BEK/2025/1/00054/DEC/1.

\end{document}